\documentclass[11pt,a4paper]{amsart}
\pdfoutput=1

\usepackage[utf8]{inputenc}

\usepackage{hyperref,graphicx,amssymb}

\DeclareMathOperator{\Aut}{Aut}
\DeclareMathOperator{\AutHS}{Aut(HS)}
\DeclareMathOperator{\HS}{HS}
\DeclareMathOperator{\Gal}{Gal}

\DeclareMathOperator{\PSL}{PSL}
\DeclareMathOperator{\PGL}{PGL}
\DeclareMathOperator{\PSU}{PSU}
\DeclareMathOperator{\Or}{O}
\DeclareMathOperator{\PSp}{PSp}

\newcommand{\QQ}{\mathbb{Q}}
\newcommand{\CC}{\mathbb{C}}

\newtheorem*{theorem}{Theorem}

\begin{document}

\title[Belyi maps having almost simple primitive monodromy groups]{Explicit Belyi maps over $\QQ$ having almost simple primitive monodromy groups}

\author{Dominik Barth}
\author{Andreas Wenz}

\address{Institute of Mathematics\\ University of Würzburg \\ Emil-Fischer-Straße 30 \\ 97074 Würzburg, Germany}
\email{dominik.barth@mathematik.uni-wuerzburg.de}
\email{andreas.wenz@mathematik.uni-wuerzburg.de}



\begin{abstract}
We present all Belyi maps $f:\mathbb{P}^1\CC\rightarrow\mathbb{P}^1\CC$ having almost simple primitive monodromy groups (not isomorphic to $A_n$ or $S_n$) containing rigid and rational generating triples of degree between 50 and 250.
This also leads to new polynomials having almost simple Galois groups over $\QQ(t)$.
\end{abstract}

\maketitle

\section{Introduction}

\noindent A triple $(x,y,z)\in S_n^3$ is called \emph{nice} if it satisfies the following conditions:
\begin{itemize}
\item $xyz = 1 $. 
\item $\left< x,y\right>$ is not a full symmetric or alternating group.
\item $\left< x,y\right>$ is primitive and almost simple.
\item $(x,y,z)$ is a genus 0 triple.
\item $(x,y,z)$ is rigid and each element is contained in a rational conjugacy class of $\langle x,y \rangle$.
\end{itemize}
Using Magma \cite{Magma} we can find all nice triples of permutation degree between $50$ and $250$.
\begin{theorem}
There are exactly $10$ nice triples $(x,y,z)$ up to simultaneous conjugation with $\left<x,y,z\right>$ being a permutation group of degree between $50$ and $250$. These triples generate the following groups:
\center{
\renewcommand{\arraystretch}{1.3}
\begin{tabular}{c|c|c}
 group        & degree & $\sharp$ nice triples  \\ \hline
$\Aut(\PSL(3,3))$   &  \phantom{0}52 & 1   \\
$\PGL(2,11)$  &  \phantom{0}55 & 2  \\
$N_{S_{56}}(\PSL(3,4))$  &  \phantom{0}56 & 1  \\
$\Aut(\PSU(3,3))$   &  \phantom{0}63 & 1  \\
$\Aut(M_{22})$     &  \phantom{0}77 & 1  \\
$\PSp(4,4){:}2$   &  \phantom{0}85 & 1  \\
$\Aut(\HS)$      & 100 & 2  \\
 $\Or^+(8,2)$       & 135 & 1  
\end{tabular}}
\end{theorem}

According to Riemann's existence theorem and the Riemann-Hurwitz formula  every nice triple $(x,y,z)$ of the previous theorem leads to a 3-branch-point covering $\tilde{f}: \mathbb{P}^1\mathbb{C}\to \mathbb{P}^1\mathbb{C}$ (also called Belyi map) with monodromy group isomorphic to $\left<x,y \right>$. Combining the rational rigidity criterion and some additional rationality considerations shows the existence of a Möbius transformation $\mu$ such that $f:=\tilde{f}\circ \mu \in \mathbb{Q}(X)$, see \cite{Malle} and \cite{voelklein}. Note that Belyi maps are unique up to inner and outer Möbius transformations. In this paper we will present Belyi maps to all nice triples of the theorem and verify their respective monodromy group.

\section{Computed Results and Verification}

We use the following definitions and facts:
Let $f = p/q \in \QQ(X)$ be of degree $n$ where $p$ and $q$ are coprime polynomials in $\QQ[X]$.
The arithmetic monodromy group of $f$ is defined as $A := \Gal(p(X)-tq(X) \in \QQ(t)[X])$, and the geometric monodromy group as $G := \Gal(p(X)-tq(X) \in \CC(t)[X])$. Both groups act transitively on the $n$ roots of $p(X)-tq(X)$ in a splitting field and we have $G \,\unlhd\, A \leq S_n$.
The subdegrees of $A$ correspond to the degrees of the irreducible factors of $p(X) - f(t)q(X) \in \QQ(t)[X]$.

In order to verify the respective geometric monodromy groups we will pursue the following idea:

The first step is to check whether the computed function $f = p/q = 1 + r/q$ is indeed a Belyi map, i.e. a three branch point covering $\mathbb{P}^1 \CC \rightarrow \mathbb{P}^1 \CC$ ramified over 0,\;1 and $\infty$. This can be done by combining the Riemann-Hurwitz genus formula with the factorization of the polynomials $p,q,r \in \QQ[X]$. We will omit this easy computation in the following.

The next step is to obtain the subdegrees of the arithmetic monodromy group $A$. In most cases we can then apply the following obvious divisibility criterion: If there is no subset of the subdegrees containing 1 adding up to a nontrivial divisor of the permutation degree, then $A$ is primitive.
If this fails we can apply Ritt's theorem which states that $A$ is imprimitive if and only if $f$ decomposes nontrivially into $f = g \circ h$ with $g,h \in \QQ(X)$. The latter can be tested by the decomposition algorithm found in \cite{af}.

By checking the Magma database for finite primitive groups having the desired subdegrees we obtain a list of possibilities for the arithmetic mono\-dromy group $A$, and for the geometric monodromy group $G$, since $G \unlhd A$. 
Fortunately, in most cases only one of these groups contains elements having the same cycle structure as $x,y$ and $z$.
This turns out to be the geometric monodromy group $G$.

Once $G$ has been determined we check if $G$ contains exactly one generating permutation triple (up to simultaneous conjugation) having the desired cycle structure. In our examples this is always the case and ensures that the ramification data over $0,1,\infty$ is indeed given by $\lbrace x,y,z\rbrace$.

Note that due to the rational rigidity of $(x,y,z)$ the splitting field of $p(X)-tq(X) \in \QQ(t)[X]$ is regular over $\QQ$, thus we have $A=G$.

\subsection*{Aut(PSL(3,3)) of degree 52}\hfill\\
We begin with the nice triple $(x,y,(xy)^{-1})$ generating $\Aut(\PSL(3,3))$ where
\begin{align*}
x\,=\;&(1, 41, 8, 9, 45, 32, 39, 44)(2, 13, 29, 21, 50, 26, 34, 6)(3, 35, 52, 30)\\
&(4, 7, 22, 18, 33, 43, 10, 38)(5, 37, 27, 42, 25, 15, 12, 24)\\
&(11, 51, 17, 47, 36, 31, 49, 40)(14, 20, 28, 48)(16, 19, 23, 46)
\end{align*}
and 
\begin{align*}
y\,=\;& (1, 20)(2, 34)(3, 7)(4, 16)(5, 17)(8, 41)(9, 13)(10, 52)(11, 40)\\
&(12, 23)(14, 29)(15, 25)(18, 33)(19, 47)(21, 35)(22, 43)(26, 42)\\
&(27, 45)(28, 36)(31, 39)(32, 49)(37, 51)(38, 48)(46, 50)
\end{align*}
of type
\renewcommand{\arraystretch}{1.5}
\begin{center}
\begin{tabular}{c|c|c|c}
& $x$ & $y$ & $(xy)^{-1}$ \\  \hline
cycle structure & $8^5.4^3$ & $2^{24}.1^4$ & $4^{10}.2^4.1^4$ \\
\end{tabular}
\end{center}
The corresponding Belyi map $f = p/q = 1 + r/q$ is given by
\begin{align*}
p(X) \,=\; & 2^2 \cdot  (4X^4 - 16X^3 - 24X^2 - 8X - 1)^ 2 \\
&\cdot (4X^4 - 16X^3 - 18X^2 - 8X - 1) \\
&\cdot (4X^4 + 8X^3 + 36X^2 + 28X + 5)^ 4 \\
&\cdot (4X^6 - 36X^5 - 24X^4 - 4X^3 + 9X^2 + 3X + 1)^ 4
\end{align*}
and
\begin{align*}
r(X) \,=\;& 3^3 \cdot (X + 1)^8 \cdot 
(2X^2 - 8X - 1)^ 8 \\ 
&\cdot (2X^2 + 1)^ 4 \cdot (6X^2 + 4X + 1)^ 8.
\end{align*}

\begin{proof}[Verification of monodromy] The factorization $p(X)-f(t)q(X) \in \QQ(t)[X]$ yields that $A$ has subdegrees $1,6,18,27$. Using the divisibility criterion we find that $A$ is primitive. As there is only one primitive permutation group of degree 52 with these subdegrees we obtain $A = \Aut(\PSL(3,3))$. Since $G$ is normal in $A$ we have $G= \PSL(3,3)$ or $G=\Aut(\PSL(3,3))$. The first case can be ruled out because $\PSL(3,3)$ does not contain elements having the same cycle structure as $(xy)^{-1}$, thus $G=A=\Aut(\PSL(3,3))$.
\end{proof}

\subsection*{PGL(2,11) of degree 55} \hfill\\
This group has exactly two nice triples $(x,y,(xy)^{-1})$ up to simultaneous conjugation. The first one is given by 
\begin{align*}
x\,=\;&(1, 33, 17, 8, 41, 32)(2, 51, 45, 36, 55, 50)(3, 43, 47, 24, 44, 25)\\
&(4, 40, 49, 52, 18, 6)(5, 7, 26, 16, 37, 23)(9, 21, 22, 38, 31, 20)\\
&(10, 48, 29, 30, 35, 19)(11, 28, 54)(12, 15, 13, 34, 42, 27)(39, 46, 53)
\end{align*}
and 
\begin{align*}
y\,=\;& (1, 26)(2, 41)(3, 51)(4, 32)(5, 29)(6, 37)(7, 20)(8, 47)(10, 44)(11, 50)\\
&(12, 52)(13, 15)(14, 35)(17, 46)(18, 23)(19, 42)(21, 48)(24, 38)\\
&(25, 45)(27, 30)(31, 39)(33, 53)(34, 36)(40, 54)(49, 55)
\end{align*}
of type
\renewcommand{\arraystretch}{1.5}
\begin{center}
\begin{tabular}{c|c|c|c}
& $x$ & $y$ & $(xy)^{-1}$ \\  \hline
cycle structure & $6^8.3^2.1^1$ & $2^{25}.1^5$ & $4^{12}.2^3.1^1$ \\
\end{tabular}
\end{center}
This triple corresponds to the Belyi map $f = p/q = 1+ r/q$ where
\begin{align*}
q(X) \,=\;& 11^4 \cdot  (2X + 1)\cdot (176X^3 + 1056X^2 + 330X + 31)^ 4 \\
& \cdot (264X^3 + 154X^2 + 22X + 1)^ 4 \cdot (352X^3 + 264X^2 + 99X + 14)^ 4 \\
&\cdot  (704X^3 + 132X^2 + 1)^ 4 \cdot  (1408X^3 + 693X^2 + 132X + 8)^ 2
\end{align*}
and 
\begin{align*}
r(X) \,=\;& -2^4 \cdot (22X^2 - 11X - 2)^ 6 \cdot (44X^2 + 22X + 3)^ 6\\
& \cdot (88X^2 + 55X + 1)^ 6 \cdot (176X^2 + 44X + 5)^ 6\\& \cdot (704X^2 + 242X + 17)^ 3.
\end{align*}

\begin{proof}[Verification of monodromy]
The polynomial $p(X)-f(t)q(X) \in \QQ(t)[X]$ factorizes into irreducible polynomials of degree $1,6,12,12,12,12$, thus $A$ must have the subdegrees $1,6,12,12,12,12$ and the primitivity of $A$ follows from the divisibility criterion. The only primitive group of degree 55 having these subdegrees is $\PGL(2,11)$. Since $G$ is normal in $\PGL(2,11)$ we have $G = \PGL(2,11)$ or $G = \PSL(2,11)$. The latter case can be ruled out because $\PSL(2,11)$ does not contain elements with the cycle structure of $y$. We find $G=A=\PGL(2,11)$.
\end{proof}
\noindent Concerning the second nice triple $(x,y,(xy)^{-1})$ in $\PGL(2,11)$ we have
\begin{align*}
x\,=\;&(2, 54, 55, 3)(4, 42, 33, 47)(5, 46, 48, 41)(6, 27, 49, 50)(7, 17, 38, 26)\\
&(8, 37, 36, 22)(9, 30, 32, 51)(10, 21, 19, 35)(11, 44, 39, 20)(12, 45, 31, 34)\\
&(13, 52, 40, 53)(14, 23, 29, 28)(15, 16, 18, 25)
\end{align*}
and
\begin{align*}
y\,=\;&(1, 22)(2, 34)(3, 41)(4, 39)(5, 50)(6, 52)(8, 42)(9, 54)(10, 27)(11, 23)\\
&(12, 24)(14, 44)(15, 49)(16, 33)(17, 35)(18, 37)(19, 32)(21, 40)(25, 46)\\
&(26, 47)(28, 48)(29, 45)(31, 38)(43, 53)(51, 55)
\end{align*}
of type
\renewcommand{\arraystretch}{1.5}
\begin{center}
\begin{tabular}{c|c|c|c}
& $x$ & $y$ & $(xy)^{-1}$ \\  \hline
cycle structure & $4^{13}.1^3$ & $2^{25}.1^5$ & $6^8.3^1.2^2$ \\
\end{tabular}
\end{center}
The computed results for this Belyi map $f= p/q = 1+ r/q$ are
\begin{align*}
q(X) \,=\;&  2^4 \cdot 
(3X - 1)^ 3 \cdot 
(2X^2 - 5X - 1)^ 6 \cdot 
(3X^2 - 2X + 4)^ 2  \\&\cdot
(12X^2 - 8X + 5)^ 6 \cdot 
(12X^4 + 6X^3 + 19X^2 - 3X + 3)^ 6
\end{align*}
and 
\begin{align*}
r(X) \,=\;& 
-3^3 \cdot 11^ 5 \cdot
(4X^3 - 4X^2 + 5X - 3)^ 4 \\& \cdot 
(6X^3 + 5X^2 + 2X + 1)^ 4 \\ & \cdot 
(12X^3 - 56X^2 + 15X - 9) \\ &\cdot 
(72X^6 - 144X^5 + 230X^4 - 134X^3 + 61X^2 - 14X + 2)^ 4.
\end{align*}

\begin{proof}[Verification of monodromy]
The subdegrees of $A$ are 1,\;4,\;6,\;8,\;12,\;24.
By applying the decomposition algorithm we see that $f$ is indecomposable and thus $A$ is primitive by Ritt's theorem.
It follows $G=\PGL(2,11)$ or $G=\PSL(2,11)$ and since there are no elements of cycle structure $2^{25}.1^5$ in $\PSL(2,11)$ we conclude again $G=A=\PGL(2,11)$.
\end{proof}

\subsection*{$\boldsymbol{N_{S_{56}}(\PSL(3,4)}$) of degree 56} \hfill\\
The only nice permutation triple $(x,y,(xy)^{-1})$ in this group is given by 
\begin{align*}
x\,=\;&(1, 36, 2, 5)(3, 47)(7, 45, 33, 22)(8, 31, 55, 14)(9, 21, 50, 48)(10, 16, 40, 39)\\
&(11, 54)(12, 19, 49, 23)(13, 41, 42, 15)(17, 56, 24, 30)(18, 53, 44, 25)\\
&(20, 52, 28, 35)(26, 29, 46, 37)(27, 34, 51, 43)
\end{align*}
and
\begin{align*}
y\,=\;&(1, 27, 36, 32)(2, 43, 50, 14, 52, 17, 30, 54)(3, 38, 47, 29, 33, 53, 16, 46)\\
&(4, 7, 22, 26, 28, 55, 44, 45)(5, 11, 24, 56, 20, 12, 41, 34)\\
&(6, 40, 18, 19, 35, 37, 10, 39)(8, 21)(9, 48, 51, 13, 15, 49, 25, 31)(23, 42)
\end{align*}
having the following cycle structure
\renewcommand{\arraystretch}{1.5}
\begin{center}
\begin{tabular}{c|c|c|c}
& $x$ & $y$ & $(xy)^{-1}$ \\  \hline
cycle structure & $4^{12}.2^2.1^4$ & $8^6.4^1.2^2$ & $2^{25}.1^6$ \\
\end{tabular}
\end{center}
The Belyi map $f= p/q = 1+ r/q$ consists of
\begin{align*}
p(X) \;=\;& (X^2 - 6X - 1)^ 2 \\ &\cdot
    (3X^4 - 468X^3 - 258X^2 - 60X - 5) \\ &\cdot
    (3X^4 + 36X^3 + 54X^2 + 60X + 19)^4 \\ &\cdot
    (3X^8 - 96X^7 - 12X^6 + 432X^5 + 1498X^4 \\&- 320X^3 - 348X^2 - 80X - 5)^4
\end{align*}
and
\begin{align*}
r(X) \;=\;& - 2^2 \cdot 5^5 \cdot (X^2 + 2X + 3)^2 \cdot (3X^2 + 6X + 1)^8\\&\cdot (3X^4 - 12X^3 + 38X^2 + 12X + 3)^8
\end{align*}

\begin{proof}[Verification of monodromy]
The subdegrees of $A$ turn out to be 1,\;10,\;45, thus $A$ is primitive. There are five finite primitive groups having these subdegrees, all between $\PSL(3,4)$ and $N_{S_{56}}(\PSL(3,4))$. As $\PSL(3,4)$ is simple these five groups are also the only possibilities for the geometric monodromy group $G \unlhd A$. However, only $N_{S_{56}}(\PSL(3,4))$ contains a generating triple with the desired cycle structure. It follows $A = G = N_{S_{56}}(\PSL(3,4))$.
\end{proof}

\subsection*{Aut(PSU(3,3)) of degree 63} \hfill\\
This group only has the nice triple $(x,y,(xy)^{-1})$ where 
\begin{align*}
x =\;
&(1, 43, 31, 39, 63, 35, 2)(3, 17, 24, 21, 20, 55, 53)(4, 29, 62, 11, 14, 45, 27)\\
&(5, 38, 23, 32, 48, 18, 51)(6, 13, 36, 47, 25, 8, 61)(7, 9, 15, 56, 34, 28, 42)\\
&(10, 33, 59, 60, 44, 19, 37)(12, 26, 52, 30, 54, 49, 41)\\&(16, 40, 57, 50, 22, 46, 58)
\end{align*}
and 
\begin{align*}
y =\;& (2, 53)(3, 9)(4, 38)(5, 29)(6, 51)(8, 25)(10, 19)(11, 20)(12, 44)(13, 43)\\
& (14, 21)(16, 40)(17, 32)(18, 39)(22, 62)(23, 30)(24, 26)(27, 59)(33, 54)
\\
& (34, 63)(35, 56)(36, 55)(41, 49)(42, 48)(45, 60)(46, 61)(47, 50)(57, 58)
\end{align*}
of type
\renewcommand{\arraystretch}{1.5}
\begin{center}
\begin{tabular}{c|c|c|c}
& $x$ & $y$ & $(xy)^{-1}$ \\  \hline
cycle structure & $7^9$ & $2^{28}.1^7$ & $4^{12}.2^6.1^3$ \\
\end{tabular}
\end{center}
The corresponding Belyi map is given by $f = p/q = 1+ r/q$ where 
\begin{align*}
r(X) =\;& - 2^8\cdot 3^{12} \cdot 
(X^2 - X + 2)^ 7 \\& \cdot 
(X^3 + 2X^2 - X - 1)^ 7 \cdot 
(X^3 + 9X^2 - X - 1)^ 7
\end{align*}
and
\begin{align*}
q(X) 
=\;& 
(X^3 + 30X^2 + 27X + 6)\\  \cdot 
&
(X^6 + 18X^5 + 93X^4 + 169X^3 + 144X^2 - 75X - 62)^2\\ \cdot 
&
(X^{12} + 15X^{11} - 15X^{10} - 332X^9 - 2766X^8 + 4002X^7 \\
& + 2002X^6 - 2496X^5 - 1215X^4 + 1047X^3 + 117X^2 - 108X + 36)^4.
\end{align*}
\begin{proof}[Verification of monodromy]
The subdegrees of $A$ are 1,\;6,\;24,\;32. Again,
by applying the decomposition algorithm we see that $A$ is primitive and only three possibilities for $A$ remain. Since $G$ is normal in $A$ there are four possibilities for $G$. Fortunately, among these groups only $\Aut(\PSU(3,3))$ contains elements of cycle structure of $y$ and $(xy)^{-1}$. This yields $G=A= \Aut(\PSU(3,3))$.
\end{proof}

\subsection*{$\boldsymbol{\Aut(M_{22})}$ of degree 77} \hfill\\
The only nice triple $(x,y,(xy)^{-1})$ is given by 
\begin{align*}
x =\;&
(1, 14, 3, 53, 31, 27, 71, 62, 10, 65, 61)(2, 50, 46, 29, 12, 7, 56, 19, 63, 28, 25)\\
&(4, 36, 38, 44, 17, 13, 66, 43, 39, 9, 72)(5, 49, 68, 51, 58, 59, 70, 15, 11, 23, 33)\\
&(6, 55, 42, 67, 32, 21, 45, 64, 48, 77, 57)(8{,}41{,}60{,}20,26, 74, 76, 24, 69, 52, 40)\\
&(16, 22, 54, 35, 34, 37, 18, 73, 75,30, 47)
\end{align*}
and
\begin{align*}
y =\;&
(1, 54)(2, 59)(3, 48)(4, 20)(6, 32)(7, 29)(11, 38)(13, 43)(14, 51)(15, 19)\\&
(18, 37)(21, 57)(22, 46)(24, 73)(30, 44)(31, 40)(33, 45)(34, 52)(35, 71)\\
&(36, 64)(39, 75)(47, 56)(49, 77)(50, 58)(53, 60)(62, 65)(63, 70)(72, 76)
\end{align*}
of type
\renewcommand{\arraystretch}{1.5}
\begin{center}
\begin{tabular}{c|c|c|c}
& $x$ & $y$ & $(xy)^{-1}$ \\  \hline
cycle structure & $11^7$ & $2^{28}.1^{21}$ & $4^{16}.2^6.1^1$ \\
\end{tabular}
\end{center}
The corresponding Belyi map is $f = p/q = 1+ r/q$ with 
\begin{align*}
p(X) =\;& 2^{22} \cdot 
(X^2 + X + 3)^{11} \\ & \cdot 
(X^5 - 3X^4 - 14X^3 + 15X^2 + X - 1)^{11}
\end{align*}
and
\begin{align*}
r(X) =\;&  - 11^4 \cdot 
(X^4 + 2X^3 + 7X^2 - 16X - 2)^ 4 \\ \cdot &
(X^6 + 14X^5 + 34X^4 + 8X^3 - 30X^2 + 60X + 16)^ 4 \\ \cdot &
(4X^6 + X^5 + 15X^4 + 10X^3 - 10X^2 - 2X - 2)^ 4 \\ \cdot &
(16X^6 - 29X^5 + 71X^4 - 136X^3 + 92X^2 - 8X - 8)^ 2.
\end{align*}

\begin{proof}[Verification of monodromy]
The subdegrees of $A$ are $1,16,60$, thus $A$ is primitive. The classification of finite primitive rank $3$ groups yields $A=M_{22}$ or $A=\Aut(M_{22})$. Since $M_{22}$ does not contain elements with the cycle structure of $y$, we see $G=A=\Aut(M_{22})$.
\end{proof}

\subsection*{PSp(4,4):2 of degree 85}\hfill \\
We now consider the nice triple $(x,y,(xy)^{-1})$ where
\begin{align*}
x =\;&
 (1, 85, 49, 26, 15, 39, 65, 24, 37, 4, 23, 3, 28, 19, 76) \\
&(2, 82, 64, 74, 52, 58, 20, 70, 43, 7, 68, 12, 53, 40, 16)\\
&(5, 81, 51, 67, 54, 44, 41, 77, 30, 21, 71, 63, 33, 66, 18)\\
&(6, 42, 46, 50, 60, 22, 73, 80, 47, 45, 14, 31, 13, 55, 79)\\
&(8, 62, 56, 36, 72, 69, 35, 25, 10, 84, 48, 34, 59, 27, 11)\\
&(9, 57, 83, 78, 17)(29, 38, 61, 75, 32)
\end{align*}
and
\begin{align*}
y =\;&
(3, 70)(5, 29)(6, 17)(7, 68)(8, 27)(9, 74)(14, 23)(15, 39)(16, 57)(18, 33)\\&
(19, 52)(20, 58)(22, 35)(24, 31)(25, 60)(26, 65)(30, 71)(32, 81)(34, 56)\\&
(36, 51)(38, 53)(40, 63)(41, 78)(42, 84)(43, 61)(44, 48)(45, 75)(46, 50)\\&
(47, 72)(49, 55)(54, 67)(64, 82)(73, 80)(76, 79)(77, 83).
\end{align*}
It has the following cycle structure:
\renewcommand{\arraystretch}{1.5}
\begin{center}
\begin{tabular}{c|c|c|c}
& $x$ & $y$ & $(xy)^{-1}$ \\  \hline
cycle structure & $15^5.5^2$ & $2^{35}.1^{15}$ & $4^{16}.2^7.1^7$ \\
\end{tabular}
\end{center}
This leads to the Belyi map $f= p/q = 1 + r/q$ with
\begin{align*}
q(X) =\;& X \cdot 
(5X^3 + 20X^2 + 20X + 6)^ 2 \\ & \cdot  
(5X^4 + 10X^3 - 14X - 10)^ 4 \\ & \cdot
(5X^4 + 10X^3 - 8X - 4)^ 2 \\& 
\cdot (5X^6 + 30X^5 + 60X^4 + 8X^3 - 48X^2 - 24X - 4) \\ & \cdot 
(625X^{12} + 3750X^{11} + 7500X^{10} + 3500X^9 - 3750X^8 - 1500X^7 
\\ & \phantom{00} + 2700X^6 + 3000X^5 + 2100X^4 + 1040X^3 + 240X^2 - 8)^ 4
\end{align*}
and
\begin{align*}
r(X) =\;& -2^ {24} \cdot 
(5X^2 + 5X + 2)^ 5 \cdot 
(5X^4 + 10X^3 - 5X - 1)^{15}.
\end{align*}

\begin{proof}[Verification of monodromy] Again, $A$ turns out to be a primitive rank $3$ group with subdegrees $1,20,64$. Thus, we either have $A=\text{PSp}(4,4)$ or $A=\text{PSp}(4,4){:}2$. As $\text{PSp}(4,4)$ does not contain elements with the cycle structure of $y$, we obtain $A=G=\text{PSp}(4,4){:}2$.
\end{proof}

\subsection*{Aut(HS) of degree 100} \hfill\\
Both nice triples $(x,y,(xy)^{-1})$ generating $\AutHS$ have already been discussed in \cite{bw}. The first one is given by
\begin{align*}
x \;=\;& (1, 23, 53, 86)(2, 36, 29, 43)(3, 15, 46, 6)(4, 80, 71, 81)(5, 75, 16, 47)
 \\& (7, 32, 60, 8)(9, 76, 100, 51)(10, 50, 49,
    34)(11, 28, 74, 84)(12, 72, 37, 52)
    \\& (13, 21, 96, 88)(14, 41, 40, 87)(17, 42, 45, 79)(18, 63, 19, 20)(22, 99, 39, 
    89)
    \\& (24, 59, 77, 38)(25, 68, 26, 35)(27, 69, 73, 48)(30, 92, 33, 82)(31, 56, 93, 58) \\&(44, 98, 67, 64)(54, 95, 85, 
    62)(55, 65, 94, 61)(57, 78, 83, 97)(66, 90, 70, 91),
\\[2mm]
y \;=\;&(1, 75, 5, 71, 15)(2, 43, 52, 89, 39)(3, 18, 100, 33, 35, 26, 58, 32, 53, 23)
\\ &(4, 81, 47, 16, 86, 7, 42, 38, 77, 
    59)(6, 41, 14, 87, 82, 76, 9, 97, 19, 63)
    \\&(8, 60, 93, 56, 13, 61, 36, 99, 70, 45)(10, 65, 55, 88, 12, 29, 94, 34,
    49, 50)
    \\& (11, 44, 64, 25, 92)(17, 72, 96, 69, 28, 30, 40, 46, 80, 24)(20, 83, 78, 57, 51)\\&(21, 31, 68, 67, 98, 84, 
    74, 27, 48, 73)(22, 37, 79, 90, 66, 95, 54, 62, 85, 91)
\end{align*}
of type
\renewcommand{\arraystretch}{1.5}
\begin{center}
\begin{tabular}{c|c|c|c}
& $x$ & $y$ & $z$ \\  \hline
cycle structure & $4^{25}$ & $10^{8}.5^4$ & $2^{35}.1^{30}$ \\
\end{tabular}
\end{center}
The Belyi map $f=p/q=1+r/q$ is given by
\begin{align*}
p(X) \;=\;
     &(7X^5 - 30X^4 + 30X^3 + 40X^2 - 95X + 50) ^4 \cdot \\ 
     &(2X^{10} - 20X^9 + 90X^8 - 240X^7 + 435X^6 - 550X^5 \\ &+ 425X^4 - 100X^3 - 175X^2 + 250X - 125)^4 \cdot \\
    & (2X^{10} + 5X^8 - 40X^6 + 50X^4 - 50X^2 + 125)^4,\\
q(X) \;=\;&
2^8 \cdot  (X^4 - 5)^5 \cdot \\
& (X^8 - 20X^6 + 60X^5 - 70X^4 + 100X^2 - 100X + 25)^{10}.    
\end{align*}

\begin{proof}[Verification of monodromy]
Once again, $A$ is a primitive rank 3 group with subdegrees $1,22,77$, therefore $A = \HS$ or $A = \AutHS$ and we find $G = \HS$ or $G = \AutHS$. Since $\HS$ does not contain elements having the cycle structure of $(xy)^{-1}$  we get 
$A = G = \AutHS$.
\end{proof}
\noindent The second triple $(x,y,(xy)^{-1})$ consists of the following permutations: 
\begin{align*}
x \;=\; & (1, 64, 8, 54, 37)(2, 20, 81, 42, 49)(3, 98, 32, 73, 89)(4, 96, 86, 15, 79)\\
   & (5, 22, 28, 78, 48)(6, 67, 97, 40, 14)(7, 58, 82, 59, 18)(9, 16, 87, 85, 60)\\
   & (10, 70, 41, 56, 55)(11, 77, 36, 25, 68)(12, 17, 19, 21, 80)(13, 35, 90, 33, 91)\\
   & (23, 50, 66, 84, 27)(24, 72, 95, 52, 76)(26, 99, 100, 57, 93)(29, 71, 38, 69, 65)\\
   & (30, 74, 94, 53, 51)(31, 45, 47, 75, 34)(43, 63, 44, 46, 62)
   \end{align*}
and
   \begin{align*}
y \;=\; & (1, 20)(2, 64)(3, 76)(4, 45)(5, 83)(6, 26)(7, 13)(8, 74)(9, 41)(10, 63)(11, 25) \\
   &(12, 66)(14, 21)(15, 52)(16, 62)(17, 33)(18, 35)(19, 42)(22, 60)(23, 58)\\    
   &(24, 73)(28, 98)(29, 82)(30, 53)(31, 61)(32, 59)(34, 67)(36, 95)(37, 85)\\
   &(38, 47)(39, 51)(40, 80)(43, 92)(44, 78)(46, 99)(48, 55)(49, 94)(50, 91)\\
   &(54, 90)(65, 88)(69, 72)(71, 75)(77, 79)(81, 87)(84, 97)(86, 100)(93, 96)
\end{align*}
of type
\renewcommand{\arraystretch}{1.5}
\begin{center}
\begin{tabular}{c|c|c|c}
& $x$ & $y$ & $z$ \\  \hline
cycle structure & $5^{19}.1^5$ & $2^{47}.1^6$ & $6^{10}.3^{10}.2^5$ \\
\end{tabular}
\end{center}
The corresponding Belyi map $f= p/q = 1 + r/q$ can be computed as
\begin{align*}
p(X) \;=\;3^3 \;\cdot\; & (X^4 - 8 X^3 - 6 X^2 + 8X + 1)^5 \cdot \\
    &(X^5 - 5X^4 + 50X^3 + 70X^2 + 25X + 3)^5 \cdot \\
    &(3X^5 - 5X^4 - 5X^3 + 35X^2 + 40X + 4) \cdot \\
    &(9X^{10} - 30X^9 + 55X^8 - 200X^7 + 210X^6 + 924X^5 \\
    &- 890X^4 - 360X^3 + 1925X^2 - 1070X + 291)^5
\end{align*}
and
\begin{align*}
q(X) \;=\;
    &(3X^5 - 35X^4 + 90X^3 - 50X^2 + 15X + 9)^2 \cdot \\
    &(9X^{10} - 120X^9 + 10X^8 - 1960X^7 - 1090X^6 + 3304X^5 \\
    &  - 760X^4 - 920X^3 + 145X^2 + 80X + 6)^3 \cdot \\
    & (3X^{10} - 10X^9 - 65X^8 + 160X^7 - 90X^6 - 932X^5 \\
    & - 330X^4 + 880X^3 + 1255X^2 + 830X + 27)^6.
\end{align*}

\begin{proof}[Verification of monodromy] One can apply the exact same proof as before to show $G=A=\AutHS$.
\end{proof}

\subsection*{$\boldsymbol{\Or^+(8,2)}$ of degree 135} \hfill\\
In this group the nice triple $(x,y,(xy)^{-1})$ is given by
\begin{align*}
x=\;
&(1, 94, 65, 71, 134, 80, 107, 98, 4)(2, 104, 58, 121, 97, 116, 88, 8, 23)\\
&(3, 69, 36, 32, 29, 73, 102, 128, 106)(5, 14, 124, 105, 67, 18, 49, 117, 34)\\
&(6, 28, 100, 41, 135, 31, 48, 109, 17)(7, 133, 112, 53, 91, 15, 25, 122, 129)\\
&(9, 62, 99, 96, 131, 77, 10, 81, 52)(11, 56, 110, 13, 115, 111, 95, 89, 54)\\
&(12, 64, 113, 108, 20, 76, 50, 22, 55)(16, 61, 83, 118, 75, 66, 39, 35, 132)\\
&(19, 85, 68, 126, 40, 125, 74, 130, 43)(21, 47, 79, 78, 72, 84, 24, 37, 57)\\
&(26, 38, 70, 90, 92, 103, 63, 120, 44)(27, 119, 127, 42, 87, 82, 101, 93, 45)\\
&(30, 59, 86, 51, 33, 60, 123, 46, 114)
\end{align*}
and
\begin{align*}
y=\; 
&(3, 118)(4, 110)(5, 132)(7, 36)(9, 33)(10, 46)(12, 112)(13, 129)(16, 65)\\
&(17, 106)(20, 113)(21, 107)(22, 55)(25, 61)(26, 27)(28, 30)(29, 37)(31, 109)\\
&(32, 98)(35, 130)(40, 42)(43, 99)(44, 125)(45, 90)(47, 49)(50, 91)(51, 93)\\
&(52, 60)(54, 56)(58, 116)(59, 128)(62, 82)(63, 73)(64, 69)(66, 96)(67, 78)\\
&(71, 117)(72, 105)(74, 124)(75, 135)(76, 83)(77, 121)(80, 134)(84, 120)\\
&(85, 87)(86, 103)(88, 104)(89, 94)(95, 122)(97, 114)(100, 131)(119, 127)
\end{align*}
and is of the following type
\renewcommand{\arraystretch}{1.5}
\begin{center}
\begin{tabular}{c|c|c|c}
& $x$ & $y$ & $(xy)^{-1}$ \\  \hline
cycle structure & $9^{15}$ & $2^{52}.1^{31}$ & $4^{30}.2^6.1^3$ \\
\end{tabular}
\end{center}
The corresponding Belyi map $f = p/q = 1  + r /q$ consists of
\begin{align*}
p(X) =\;& 2^{22} \cdot  (3X^3 - 9X^2 - 9X - 2)^9 \cdot (3X^3 + 9X^2 + 6X + 1)^9 \\
&\cdot (27X^9 + 243X^8 + 567X^7 + 513X^6 + 162X^5 \\&- 27X^4 + 9X^3 + 27X^2 + 9X + 1)^9
\end{align*}
and
\begin{align*}
q(X) =\;& (3X^3 - 9X - 2) \cdot (3X^3 + 27X^2 + 27X + 7)^2 \cdot (6X^3 + 9X^2 - 1)^2\\
&\cdot (36X^6 + 189X^5 + 189X^4 + 96X^3 + 36X^2 + 9X + 1)^4\\
&\cdot(81X^{12} + 1944X^{11} + 11178X^{10} + 27648X^9 + 29403X^8- 1944X^7\\ 
& - 39150X^6 - 44712X^5 - 25434X^4 - 8088X^3 - 1332X^2 - 72X + 4)^4\\
&\cdot (648X^{12} + 3888X^{11} + 11907X^{10} + 15120X^9 + 13365X^8+ 14580X^7\\ 
& + 11772X^6 + 3240X^5 - 1782X^4 - 1632X^3 - 504X^2 - 72X - 4)^4.
\end{align*}

\begin{proof}[Verification of monodromy] We find that $A$ is a primitive rank 3 group with subdegrees $1,64,70$. Therefore we have $A = \Or^+(8,2)$ or $A=\Or^+(8,2).2$ and due to normality $G = \Or^+(8,2)$ or $G=\Or^+(8,2).2$. Note that $G$ is generated by permutations of cycle structure $9^{15}$, $2^{52}.1^{31}$, $4^{30}.2^6.1^3$ and by inspecting the sizes of conjugacy classes of $\Or^+(8,2)$ and $\Or^+(8,2).2$ we can conclude that there are no elements with these cycle structures in $\Or^+(8,2).2 \setminus \Or^+(8,2)$. It follows $G = \Or^+(8,2)$.
Because $G$ contains only one generating triple (up to simultaneous conjugation) having the desired cycle structure the ramification data over $0,1,\infty$ is given by $\lbrace x,y,(xy)^{-1} \rbrace$. Since $(x,y,(xy)^{-1})$ is rigid we have $A=G=\Or^+(8,2)$.
\end{proof}

\section*{Acknowledgements}
We would like to thank Peter Müller for introducing us to the subject of this work and for providing us with helpful suggestions.

\end{document}